\documentclass[12pt]{amsart}

\usepackage[dvips]{graphics}

\begin{document}

\title[Compact Minimal Surfaces without Area Bounds]{Compact Embedded Minimal Surfaces of Positive Genus without Area Bounds}
\author{Brian Dean}
\address{Department of Mathematics\\
    Johns Hopkins University\\
    3400 N. Charles St.\\
    Baltimore, MD  21218}
\thanks{Johns Hopkins University; bdean@math.jhu.edu; Running headline:  Compact Minimal Surfaces without Area Bounds; AMS Subjects:  53A10, 53C42; Key words:  minimal surfaces, stability, differential geometry.}
\date{}
\maketitle

\begin{abstract}
Let $M^3$ be a three-manifold (possibly with boundary).  We will
show that, for any positive integer $\gamma$, there exists an open
nonempty set of metrics on $M$ (in the $C^2$-topology on the space
of metrics on $M$) for each of which there are compact embedded
stable minimal surfaces of genus $\gamma$ with arbitrarily large
area.  This extends a result of Colding and Minicozzi, who proved
the case $\gamma =1$.
\end{abstract}

\section*{Introduction}\label{sec:int}

Throughout this paper, we use the $C^2$-topology on the space of metrics on a manifold.  Our main result is the following theorem.

\medskip
\noindent
\textbf{Theorem 1:}  Let $M^3$ be a three-manifold (possibly with boundary), and let $\gamma$ be a positive integer.  There exists an open nonempty set of metrics on $M$ for each of which there are compact embedded minimal surfaces of genus $\gamma$ with arbitrarily large area.  In fact, these can be chosen to be stable, i.e., with Morse index zero.
\medskip

Although the theorem ensures that there are ``many'' metrics for which we can embed compact genus $\gamma$ minimal surfaces of arbitrarily large area, the result is false for a large class of metrics.  Namely, a result of Choi and Wang (see [CW]) asserts that for any metric in which $M$ has Ricci curvature bounded below by a positive constant, there is an upper bound on the area of compact embedded minimal surfaces of genus $\gamma$, depending on $\gamma$ and the lower bound for $\mbox{Ric}_M$.

Colding and Minicozzi (see~\cite{CM1}) have already proved Theorem 1 for $\gamma =1$.  In Section~\ref{sec:genus2}, we will prove the theorem for $\gamma =2$, with an argument borrowing heavily from the genus one case.  The theorem will then be extended easily to genus greater than two in Section~\ref{sec:highergenus}.

It remains an open question as to whether or not the theorem remains valid for genus zero, i.e., embedded minimal 2-spheres.

\section{The genus 2 case}\label{sec:genus2}

\begin{figure}[!h!t]
   \begin{center}
      \scalebox{0.3}{\includegraphics{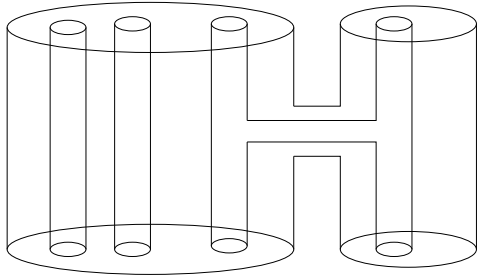}}
   \end{center}
   \caption{$\Omega_2$.  The top and bottom of the picture are identified.}
   \label{figomega2}
\end{figure}

Let $\Sigma_2$ denote the standard genus two surface.  This has fundamental group

\[
\pi_1(\Sigma_2)=\,<x_1,y_1,x_2,y_2\,|\,x_1y_1x_1^{-1}y_1^{-1}x_2y_2x_2^{-1}y_2^{-1}>
\]

\noindent where $x_1$ and $x_2$ are freely homotopic to meridians
of the two handles, where the meridians have the same orientation,
and $y_1$ and $y_2$ are freely homotopic to lines of latitude of
the two handles, where the lines of latitude have the same
orientation. Let $\Omega_2$ be a solid genus two surface with a
solid genus two surface and two solid tori removed, where the
solid tori lie in the same handle of the ambient genus two
surface.  This can be pictured as in Figure~\ref{figomega2}, where
the top and bottom of the picture are identified.

The fundamental group of $\Omega_2$ is
\[
\pi_1(\Omega_2)=\,<a,b,c_1,d_1,c_2,d_2\,|\,ad_1a^{-1}d_1^{-1}, \,bd_1b^{-1}d_1^{-1}, \,c_1d_1c_1^{-1}d_1^{-1}c_2d_2c_2^{-1}d_2^{-1}>
\]

\noindent
where the generators are as follows:

\begin{description}
\item[(i)] $a$ and $b$ are freely homotopic to meridians of the
two removed solid tori (clockwise rotation around the two removed
solid tori in Figure 1). \item[(ii)] $c_1$ is freely homotopic to
a meridian of the handle of the removed solid genus two in the
same handle of the ambient solid genus two as the removed solid
tori (clockwise rotation around the left handle of the removed
solid genus two in Figure 1). \item[(iii)] $d_1$ is freely
homotopic to a line of latitude of the left handle of the ambient
solid genus two in Figure 1. \item[(iv)] $c_2$ is freely homotopic
to a meridian of the removed solid genus two in the other handle
from the meridian which is freely homotopic to $c_1$ (clockwise
rotation around the right handle of the removed solid genus two in
Figure 1). \item[(v)] $d_2$ is freely homotopic to a line of
latitude of the right handle of the ambient solid genus two in
Figure 1, with the same orientation as the line of latitude which
is freely homotopic to $d_1$.
\end{description}

\medskip
Before we give the proof of Theorem 1 for the genus two case, we need the following proposition, whose proof is inspired by a calculation in~\cite{Es}.

\medskip
\noindent
\textbf{Proposition 1:}  Let $\Omega^n$ be a compact Riemannian manifold with boundary and dimension $n\ge 3$.  Then, the set of metrics on $\Omega$ in which $\Omega$ is strictly mean convex is open and nonempty.

\medskip
\noindent
\textbf{Proof:}  The set of such metrics is clearly open, by the definition of strictly mean convex.  To show it is nonempty, let $g$ be any metric on $\Omega$, and let $\widetilde g = e^{2f}g$ be a metric conformally related to $g$.  Let $\{e_1,...,e_n\}$ be a framing for $\Omega$ so that $g_{ij}=\delta_{ij}$ and $e_n$ is the unit normal to $\partial\Omega$ in $g$ (and therefore, $e^{-f} e_n$ is the unit normal to $\partial\Omega$ in $\widetilde g$).  Fix a point $p\in\partial\Omega$, and choose coordinates $\{x_1,...,x_n\}$ at $p$ so that, at $p$, $e_i = \frac{\partial}{\partial x_i}$ for all $i$.  Then, the second fundamental form of $\partial\Omega$ in $g$ at $p$ is given by, for $i,j=1,...,n-1$,
\[
h_{ij} = g(\nabla_{e_i}e_j,e_n) = \sum_{k=1}^{n}g(\Gamma_{ij}^k e_k,e_n) = \Gamma_{ij}^n,
\]
and the mean curvature of $\partial\Omega$ in $g$ at $p$ is given by
\[
h = \frac{1}{n-1}\sum_{i,j=1}^{n-1}g^{ij}h_{ij} = \frac{1}{n-1}\sum_{i,j=1}^{n-1}\delta_{ij}\Gamma_{ij}^n = \frac{1}{n-1}\sum_{i=1}^{n-1}\Gamma_{ii}^n
\]

The second fundamental form of $\partial\Omega$ in $\widetilde g$ at $p$ is given by, for $i,j=1,...,n-1$,
\[
\widetilde h_{ij} = \widetilde g(\widetilde\nabla_{e_i}e_j,e^{-f}e_n) = e^{2f}g(\widetilde\nabla_{e_i}e_j,e^{-f}e_n) = \sum_{k=1}^n e^f g(\widetilde\Gamma_{ij}^k e_k,e_n) = e^f\widetilde\Gamma_{ij}^n.
\]
Now,

\[
\widetilde\Gamma_{ij}^n=\frac{1}{2}\sum_{l=1}^n(\widetilde g_{jl,i}+\widetilde g_{il,j}-\widetilde g_{ij,l})\,\widetilde g^{ln}
\]

\noindent
where, for example,

\[
\widetilde g_{jl,i}=\frac{\partial}{\partial x_i}\,\widetilde g_{jl}
\]

\noindent
So, at $p$,
\begin{eqnarray*}
\widetilde\Gamma_{ij}^n &=& \frac{1}{2}e^{-2f}\sum_{l=1}^n(\widetilde g_{jl,i}+\widetilde g_{il,j}-\widetilde g_{ij,l})\,g^{ln} \\
   &=& \frac{1}{2}e^{-2f}\sum_{l=1}^n(\widetilde g_{jl,i}+\widetilde g_{il,j}-\widetilde g_{ij,l})\,\delta_{ln} \\
   &=& \frac{1}{2}e^{-2f}(\widetilde g_{jn,i}+\widetilde g_{in,j}-\widetilde g_{ij,n}) \\
   &=& \frac{1}{2}e^{-2f}\left(2\frac{\partial f}{\partial x_i}e^{2f}g_{jn}+e^{2f}g_{jn,i}+2\frac{\partial f}{\partial x_j}e^{2f}g_{in}+e^{2f}g_{in,j}\right. \\
   & & \left. {}-2\frac{\partial f}{\partial x_n}e^{2f}g_{ij}-e^{2f}g_{ij,n}\right) \\
   &=& \frac{1}{2}(g_{jn,i}+g_{in,j}-g_{ij,n})+\frac{\partial f}{\partial x_i}\delta_{jn}+\frac{\partial f}{\partial x_j}\delta_{in}-\frac{\partial f}{\partial n}\delta_{ij} \\
   &=& \Gamma_{ij}^n-\frac{\partial f}{\partial n}\delta_{ij}
\end{eqnarray*}
since $i,j<n$, where $\frac{\partial f}{\partial n}$ is the normal derivative of $f$ with respect to the unit normal $e_n$.  Therefore, we have
\[
\widetilde h_{ij} = e^f\widetilde\Gamma_{ij}^n = e^f\Gamma_{ij}^n-e^f\frac{\partial f}{\partial n}\delta_{ij} = e^f\Gamma_{ij}^n-\frac{\partial}{\partial n}(e^f)\delta_{ij}.
\]

The mean curvature of $\partial\Omega$ in $\widetilde g$ at $p$ is then given by
\begin{eqnarray*}
\widetilde h &=& \frac{1}{n-1}\sum_{i,j=1}^{n-1}\widetilde g^{ij}\widetilde h_{ij} \\
   &=& \frac{1}{n-1}\sum_{i,j=1}^{n-1}e^{-2f}\delta_{ij}\left(e^f\Gamma_{ij}^n-\frac{\partial}{\partial n}(e^f)\delta_{ij}\right) \\
   &=& \frac{1}{n-1}\sum_{i=1}^{n-1}\left(e^{-f}\Gamma_{ii}^n-e^{-2f}\frac{\partial}{\partial n}(e^f)\right) \\
   &=& e^{-f}\left(\frac{1}{n-1}\sum_{i=1}^{n-1}\Gamma_{ii}^n\right)-e^{-2f}e^f\frac{\partial f}{\partial n} \\
   &=& e^{-f}\left(h-\frac{\partial f}{\partial n}\right)
\end{eqnarray*}
We have shown that this relation holds at an arbitrarily chosen point of $\partial\Omega$, and so it holds everywhere on $\partial\Omega$ since all quantities involved are tensorial.  Let $m$ be the minimum of $h$ on $\partial\Omega$, which exists since $\partial\Omega$ is compact.  Choose $f$ so that $\frac{\partial f}{\partial n}<m$ everywhere on $\partial\Omega$ and $f\equiv 0$ outside a small tubular neighborhood around $\partial\Omega$.  Then, $\widetilde g=g$ except for a small tubular neighborhood around $\partial\Omega$, and $\widetilde h>0$ everywhere on $\partial\Omega$, so $\widetilde g$ is a metric in which $\Omega$ is strictly mean convex.  This completes the proof of Proposition 1.

\bigskip
To prove Theorem 1, we will also need the following lemma.

\medskip
\noindent
\textbf{Lemma 1:}  Let $N^3$ be a compact Riemannian manifold, and let $\{M_n\}\subset N$ be a sequence of stable, compact, connected, embedded minimal surfaces without boundary such that the following conditions hold:
\begin{description}
\item[(i)] there exists a constant $C_1>0$ such that Area$(M_n)\le C_1$ for all $n$.
\item[(ii)] there exists a constant $C_2>0$ such that

\[
\sup_{M_n}|A_n|^2\le C_2
\]
for all $n$, where $A_n$ is the second fundamental form of $M_n$.
\end{description}
Then, a subsequence of $\{M_n\}$ converges to a compact, connected, embedded minimal surface without boundary $M\subset N$ of finite multiplicity.

\medskip
\noindent
\textbf{Proof:}  Take a finite covering $\{B_r(y_j)\}$ of $N$ so that $\{B_{r/2}(y_j)\}$ is still a covering of $N$.  Then, by~\cite{CM2}, a subsequence of $\{M_n\}$ converges in each $B_{r/2}(y_j)$ to a lamination $M$ with minimal leaves.  By taking a diagonal subsequence, we have a subsequence of $\{M_n\}$, which we still call $\{M_n\}$, converging to $M$ everywhere.  $M$ is clearly minimal, and it is embedded by the maximum principle.

We claim that the number of leaves of $M_n$ in each $B_r(y_j)$ which intersect $B_{r/2}(y_j)$ has an upper bound which is uniform in $n$ and $j$.  Let $\Gamma_{n,j}$ be any such leaf.  Then, there exists $x_j\in\Gamma_{n,j}\cap\partial B_{r/2}(y_j)$.  So, $B_{r/2}(x_j)\subset B_r(y_j)$, and by monotonicity of area, there exists a constant $C>0$ so that

\[
\mbox{Area}(\Gamma_{n,j}\cap B_{r/2}(x_j))\ge C\left(\frac{r}{2}\right)^2.
\]
  So, each $\Gamma_{n,j}$ is of at least some fixed positive area, and so the area bound $(i)$ gives an upper bound for the number of such leaves which is uniform in $n$ and $j$.  We can take a subsequence so that the number of leaves of $M_n$ is the same in each $B_{r/2}(y_j)$ for all $n,j$.  Then, the limit $M$ must have finite multiplicity, although the multiplicity may be different in each connected component of $M$.  We have shown that each connected component of $M$ is a closed surface.  The diameter of $M$ is bounded, since $M$ is covered by finitely many balls $B_{r/2}(y_j)\cap M$.  So, $M$ is compact.  $M$ is without boundary since each $M_n$ is without boundary.

It remains to show that $M$ is connected, which would imply that $M_n\rightarrow M$ with fixed finite multiplicity.  Suppose $M$ is not connected, and let $A$ and $B$ be distinct connected components of $M$.  Then, $\epsilon=\mbox{dist$(A,B)$}>0$.  Let $R=\left\{x\in N|\,\frac{\epsilon}{3}<\mbox{dist$(x,A)$}<\frac{2\epsilon}{3}\right\}$.  So, $R$ is disjoint from both $A$ and $B$.  Since $M_n\rightarrow M$, for large enough $n$ we have $M_n\cap A\ne\emptyset$ and $M_n\cap B\ne\emptyset$, but $M_n\cap R=\emptyset$, contradicting the connectedness of $M_n$.  So, $M$ is connected.

Therefore, $M_n$ converges to a compact, connected, embedded minimal surface without boundary $M\subset N$ of finite multiplicity.  This completes the proof of Lemma 1.

\bigskip
\noindent
\textbf{Proof of Theorem 1 for $\gamma =2$:}

Given any three-manifold $M$, we can embed $\Omega_2$ in $M$.  Choose a metric $g$ on $M$ so that $\Omega_2$ is strictly mean convex (by Proposition 1, the set of such $g$'s is open and nonempty).  Let $f_n:\Sigma_2\rightarrow\Omega_2$ be a map such that the induced map ${f_n}_\#:\pi_1(\Sigma_2)\rightarrow\pi_1(\Omega_2)$ is the following:
\begin{eqnarray*}
{f_n}_\#(x_1) &=& (ba)^nc_1(ba)^{-n}a^{-1}ba \\
{f_n}_\#(y_1) &=& d_1 \\
{f_n}_\#(x_2) &=& c_2 \\
{f_n}_\#(y_2) &=& d_2
\end{eqnarray*}
It is easy to see that there exist such maps $f_n$ which are embeddings.

One can check that ${f_n}_\#(x_1)$ minimizes the word metric for its conjugacy class:  by conjugating ${f_n}_\#(x_1)$ by any element of $\pi_1(\Omega_2)$ and using the relations of $\pi_1(\Omega_2)$, one can not decrease the length of ${f_n}_\#(x_1)$ in the word metric (for the definition of word metric, see~\cite{CM1}).

So, for each $n$, we have an embedded incompressible genus two surface $\Sigma_{2,n}=f_n(\Sigma_2)$.  By~\cite{ScY}, there are immersed least-area (minimal) genus two surfaces $\Gamma_{2,n}\subset\Omega_2$ with $\Gamma_{2,n}\cap\partial\Omega_2=\emptyset$ so that $\Gamma_{2,n}$ and $\Sigma_{2,n}$ induce the same mapping from $\pi_1(\Sigma_2)$ to $\pi_1(\Omega_2)$ for each $n$.  Since the $\Sigma_{2,n}$ are embedded,~\cite{FHS} implies that the $\Gamma_{2,n}$ are embedded.

We claim that the areas of the $\Gamma_{2,n}$'s are unbounded.  Assume not.  Then, there exists a constant $C_1>0$ such that Area$(\Gamma_{2,n})\le C_1$ for all $n$.  The $\Gamma_{2,n}$ are stable since they are area-minimizing.  So, by~\cite{Sc}, we get a uniform curvature estimate:  there exists a constant $C_2>0$ such that, for small enough $r$ and all $\sigma\in(0,r]$,

\[
\sup_{B_{r-\sigma}}|A_n|^2\le\frac{C_2}{\sigma^2}
\]
for all $n$ and all balls $B_{r-\sigma}\subset\Omega_2$, where $A_n$ is the second fundamental form of $\Gamma_{2,n}$.  Since the $\Gamma_{2,n}$ are all without boundary, we get a uniform curvature estimate on all of $\Gamma_{2,n}$, instead of just on balls.  Therefore, by Lemma 1, a subsequence of $\{\Gamma_{2,n}\}$ converges to a compact, connected, embedded minimal surface without boundary $\Gamma_2\subset\Omega_2$ of finite multiplicity.

For large $n$, the $\Gamma_{2,n}$ are coverings of $\Gamma_2$ by the maximum principle, and the degree of the covering is proportional to $n$.  Let $n\rightarrow\infty$.  Then, $\Gamma_2$ has infinite multiplicity, a contradiction.  Therefore, the areas of the $\Gamma_{2,n}$'s are unbounded.  This completes the proof of Theorem 1 for the case $\gamma =2$.

\section{The general case: $\gamma\ge 2$}\label{sec:highergenus}

We now move to the general case.  The arguments for fixed genus $\gamma\ge 2$ are essentially the same as in the genus 2 case.

\begin{figure}[!h!t]
   \begin{center}
      \scalebox{0.3}{\includegraphics{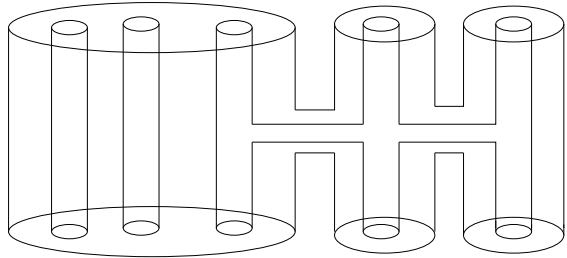}}
   \end{center}
   \caption{$\Omega_3$.  The top and bottom of the picture are identified.}
   \label{figomega3}
\end{figure}

Let $\Sigma_{\gamma}$ denote the standard genus $\gamma$ surface, $\gamma\ge 2$.  This has fundamental group

\[
\pi_1(\Sigma_{\gamma})=\,<x_1,y_1,\ldots,x_{\gamma},y_{\gamma}\,|\,[x_1y_1]\cdots[x_{\gamma}y_{\gamma}]>
\]

\noindent where the $x_i$ are freely homotopic to meridians of the
handles, all with the same orientation, the $y_i$ are freely
homotopic to lines of latitude of the handles, all with the same
orientation, and $[x_iy_i]=x_iy_ix_i^{-1}y_i^{-1}$ for
$i=1,\ldots,\gamma$. Let $\Omega_{\gamma}$ be a solid genus
$\gamma$ surface with a solid genus $\gamma$ surface and two solid
tori removed, where the solid tori both lie in one of the end
handles of the ambient genus $\gamma$ surface (the case $\gamma
=3$ is shown in Figure~\ref{figomega3}, where the top and bottom
of the picture are identified).

The fundamental group of $\Omega_{\gamma}$ is

\[
\pi_1(\Omega_{\gamma})=\,<a,b,c_1,d_1,\ldots,c_{\gamma},d_{\gamma}\,|\,[ad_1],\,[bd_1],\,[c_1d_1]\cdots[c_{\gamma}d_{\gamma}]>
\]

\noindent where the generators are defined as in the case $\gamma
=2$ (so, $a$, $b$, and all $c_i$ are freely homotopic to meridians
with the same orientation, and all $d_i$ are freely homotopic to
lines of latitude with the same orientation).

\medskip
\noindent
\textbf{Proof of Theorem 1:}

Given any three-manifold $M$, we can embed $\Omega_{\gamma}$ in $M$.  Choose a metric $g$ on $M$ so that $\Omega_{\gamma}$ is strictly mean convex (by Proposition 1, the set of such $g$'s is open and nonempty).  Let $f_n:\Sigma_{\gamma}\rightarrow\Omega_{\gamma}$ be a map such that the induced map ${f_n}_\#:\pi_1(\Sigma_{\gamma})\rightarrow\pi_1(\Omega_{\gamma})$ is the following:
\begin{eqnarray*}
{f_n}_\#(x_1) &=& (ba)^nc_1(ba)^{-n}a^{-1}ba \\
{f_n}_\#(x_i) &=& c_i \mbox{ for $i=2,\ldots,\gamma$} \\
{f_n}_\#(y_i) &=& d_i \mbox{ for $i=1,\ldots,\gamma$}
\end{eqnarray*}
It is easy to see that there exist such maps $f_n$ which are embeddings.

The proof then proceeds exactly as in the genus 2 case. The results of~\cite{ScY},~\cite{FHS}, and~\cite{Sc} again apply.

\end{document}